\documentclass[12pt,a4paper]{amsart}
\usepackage[top=35mm, bottom=30mm, left=30mm, right=30mm]{geometry}
\usepackage{mathptmx}
\usepackage{mathrsfs}
\usepackage{verbatim}
\usepackage{url}
\usepackage[all]{xy}
\usepackage{xcolor}
\usepackage[colorlinks=true,citecolor=blue]{hyperref}
\usepackage{amsmath,amsthm}
\usepackage{amssymb}
\usepackage{enumerate}
\usepackage{graphicx}

\makeatletter
\@namedef{subjclassname@2020}{%
    \textup{2020} Mathematics Subject Classification}

\theoremstyle{definition}

\numberwithin{equation}{section}
\begin{document}


\baselineskip=17pt


\title{Orientation Preserving Homeomorphisms of the Plane having BP-Chain Recurrent Points}

\author[J.~Mai]{Jiehua Mai}
\address{School of Mathematics and Quantitative
Economics $\&$ Guangxi Key Laboratory of Quantitative Economics,
Guangxi University of Finance and Economics, Nanning, Guangxi,
530003, P. R. China; Institute of Mathematics, Shantou University,
Shantou, Guangdong, 515063, P. R. China} \email{jiehuamai@163.com}

\author[K.~Yan]{Kesong Yan*}
\thanks{*Corresponding author}
\address{School of Mathematics and Statistics, Hainan Normal
University,  Haikou, Hainan, 571158, P. R. China}
\email{ksyan@mail.ustc.edu.cn}

\author[F.~Zeng]{Fanping Zeng}
\address{School of Mathematics and Quantitative
Economics $\&$ Guangxi Key Laboratory of Quantitative Economics,
Guangxi University of Finance and Economics, Nanning, Guangxi,
530003, P. R. China} \email{fpzeng@gxu.edu.cn}

\begin{abstract}
More than a century ago, L. E. J. Brouwer proved a famous theorem,
which says that any orientation preserving homeomorphism of the
plane having a periodic point must have a fixed point. In recent
years, there are still some authors giving various proofs of this
fixed point theorem. In \cite{Fa}, Fathi showed that the condition
``having a periodic point'' in this theorem can be weakened to
``having a non-wandering point''. In this paper, we first give a new
proof of Brouwer's theorem, which is relatively more simpler and the
statement is more compact. Further, we propose a notion of BP-chain
recurrent points, which is a generalization of the concept of
non-wandering points, and we prove that if an orientation preserving
homeomorphism of the plane has a BP-chain recurrent point, then it
has a fixed point. This further weakens the condition in the
Brouwer's fixed point theorem on plane.
\end{abstract}


\maketitle

\pagestyle{myheadings} \markboth{J. H. Mai, K. S. Yan, F. P. Zeng
}{Orientation Preserving Homeomorphisms of the Plane having BP-Chain
Recurrent Points}

\section{Introduction}

More than a century ago, L. E. J. Brouwer \cite{Bro12} first stated
and proved the following famous theorem:

\vspace{3mm}{\textbf{Theorem A.}} \ {\it Any continuous map from the
closed unit ball in $n$-dimensional Euclidean space to itself has a
fixed point.}

\vspace{3mm}

For more than a century, the attention to Theorem A has never been
interrupted. There are various different proofs, generalizations and
applications of Theorem A, see the survey papers \cite{Maw07, Maw20,
Park99, Sub18} and references therein. Until recently, there are
still some authors giving new and simple proofs of Theorem A and its
generalizations, for example, see \cite{BM22, IKM14, IKM21, Maw19,
Su97, Vol08}.

\vspace{3mm}

In \cite{Bro12b}, Brouwer proved another famous fixed point theorem:

\vspace{3mm}{\textbf{Theorem B.}} \ {\it Any orientation preserving
homeomorphism of the plane having a periodic point must have a fixed
point.}

\vspace{3mm}

Similar to Theorem A, Theorem B has gained continuous attention by
many people, and some new and ingenious proofs have been given in
the last few years. A simple modern proof due to M. Brown can be
found in \cite{B1} and A. Fathi \cite{Fa} has given another very
nice proof. In \cite{B1}, Brown proved Theorem B by using the
Alexander isotopy, which relied on the notion of index and the
notion of free equivalence. In \cite{Fa}, Fathi showed that the
condition ``having a periodic point'' in Theorem B can be weakened
to ``having a non-wandering point''. In \cite{F2}, Franks
investigated fixed point free orientation preserving homeomorphisms
of the plane and gave a short new proof of Brouwer plane translation
theorem. In addition, Barge and Franks \cite{BF} further discussed
fixed point free orientation preserving homeomorphisms of the plane
and the recurrence of connected subsets. They gave a short
elementary proof of Theorem B which is similar to the proof of Fathi
\cite{Fa}. In \cite{B2}, Brown showed that Theorem B follows from a
particularly simple and elegant application of the notion of index
of a homeomorphism along an arc.

\vspace{3mm}

In this paper, we first give a new proof of Theorem B, which is
relatively more simpler and the statement is more compact. Further,
we propose a notion of BP-chain recurrent points, which is a
generalization of the concept of non-wandering points. After giving
an example to show that the condition ``having a periodic point'' in
Theorem B can not be weaken to ``having a  chain recurrent point'',
we prove that if an orientation preserving homeomorphism of the
plane has a BP-chain recurrent point then it has a fixed point. This
further weakens the condition in Theorem B.

\section{A simple proof of Brouwer's fixed
point theorem for plane homeomorphisms}

Let\;\! $\mathbb{N}$ \;\!be the set of all positive integers, \,and
let\;\! $\mathbb{Z}_+=\mathbb{N}\cup\{0\}$. For any
\;\!$n\in\mathbb{N}$\;\!, \,write \;\!$\mathbb{N}_n=\{1,\cdots ,
n\}$. For any $x \in \mathbb{R}^2$ and $r>0$, write $B^2(x,r)=\{y
\in \mathbb{R}^2: d(x,y) \leq r\}$, $S^1(x,r)=\{y \in \mathbb{R}^2:
d(x,y)=r\}$, and $\mathrm{Int}\,B^2(x,r)=\{y \in \mathbb{R}^2:
d(x,y)<r\}$. For any $V \subset \mathbb{R}^2$, write
$B^2(V,r)=\bigcup\{B^2(x,r): x \in V\}$. For any $X \subset
\mathbb{R}^2$, denote by $\mathrm{diam}\, X$ the diameter of $X$. If
there is a homeomorphism $h: B^2(0,1) \rightarrow X$, then $X$ is
called a {\it\textbf{disk}}, and we write
$\stackrel{\circ}X\;=\mathrm{Int}\,X=h\left(\mathrm{Int}\,B^2(0,1)\right)$,
$\partial X=h\left(S^1(0,1)\right)$. If there is a homeomorphism $h:
S^1(0,1) \rightarrow X$, then $X$ is called a {\it\textbf{circle}}.
If there is a homeomorphism $h: [0,1] \rightarrow X$, then $X$ is
called an {\it\textbf{arc}}, and we write
$\stackrel{\circ}X\;=\mathrm{Int}\,X=h\left((0,1)\right)$, $\partial
X=h\left(\{0,1\}\right)$. For any $x, y \in \mathbb{R}^2$ with $x
\neq y$, let $[x,y]=[y,x]$ be the straight line segment with
endpoints $x$ and $y$, and let $[x,x]=\{x\}$. For any arc $A \subset
\mathbb{R}^2$ and any $x, y \in A$, denote by $[x,y]_A=A[x,y]$ the
least connected subset of $A$ containing $\{x,y\}$, and write
$(x,y]_A=[y, x)_A=[x,y]_A-\{x\}$, $(x,y)_A=A(x,y)=(x,y]_A-\{y\}$.

\vspace{3mm}For any metric space $(X,d)$\;\!, denote by
\;\!$C^{\;\!0}(X)$\;\! the set of all continuous maps from $X$ to
$X$, and let $id=id_X$ be the identity map of $X$. For any $f \in
C^0(X)$, let $f^0=id$, and let $f^{\;\!n}=f\circ f^{\;\!n-1}$ be the
composition map of $f$ and $f^{\;\!n-1}$. A point $x\in X$ is called
a {\it\textbf{fixed point}} \,of $f$ if \;\!$f(x)=x$\;\!. Denote by
$\mathrm{Fix}(f)$ the set of fixed points of $f$. We say that\,$x\in
X$ is called a {\it\textbf{periodic point}}\, of $f$ if
$f^{\,n}(x)=x$ \,for some $n\in\mathbb{N}$\;\!; the smallest such
$n$ is called the {\it\textbf{period}}\, of \,$x$ under $f$. A
periodic point of period $n$ is also called an
{\it\textbf{$n$-periodic point}}.

\vspace{3mm} In \cite{BF, B1, B2, Fa, F1, F2, GK}, the authors give
variant proofs of Theorem B. A way of the proof of Theorem B is to
prove the following two propositions\,:

\vspace{3mm}\textbf{Proposition C.} \ {\it Let $f:\Bbb R^2\to
\mathbb{R}^2$ be an orientation preserving homeomorphism. If $f$ has
a $2$-periodic point,\, then $f$ has a fixed point}\;\!.

\vspace{3mm}\textbf{Proposition D.} \ {\it Let $f:\Bbb R^2\to\Bbb
R^2$ be an orientation preserving homeomorphism. If $f$ has an\,
$n$-periodic point\;\! $x$,\,$n\ge3$\;\!, \,then there exists an
orientation preserving homeomorphism \,$g:\Bbb R^2\to\Bbb R^2$ such
that\; ${\rm Fix}(g)=\,{\rm Fix}(f)$\;\!, and\, $g$ \;\!has a\;\!
$k$-periodic point for some\, $k\in \mathbb{N}_{n-1}$}\;\!.

\vspace{3mm} In the proof of Theorem 1.1 of \cite{BF},\;\! the
authors actually gave a proof of Proposition D, see
\cite[pp.\,188\;--\;189]{BF}. Relatively, the proof of Proposition D
is easier than that of Proposition C. Thus the proof of Proposition
D is a key of that of Theorem B. In \cite{BF}, \cite{B2}, \cite{Fa}
and \cite{GK}\;\!, the authors gave variant proofs of Proposition D,
respectively.

\vspace{3mm} In this section we will give shorter proofs of
Propositions C and D. This means that we give a new proof of Theorem
B. First we propose the following lemma.

\vspace{3mm}\textbf{Lemma 2.1.} \ {\it Let $f:\Bbb R^2\to\Bbb R^2$
be an orientation preserving homeomorphism. If there exists an arc
$A$ in $\mathbb{R}^2$ such that $\partial A=\{x,y\}$, $f(x)=y$,
$f(y)=x$, and $f\left(\mathrm{Int}\,A\right) \cap
\mathrm{Int}\,A=\emptyset$, then $f$ has a fixed point.}

\vspace{3mm}\textbf{Proof.} \ Write $A'=f(A)$, and $C=A' \cup A$.
Then $C$ is a circle. For any circle $Q$ in $\mathbb{R}^2$, let
$\mathrm{Dsc}(Q)$ be the disc in $\mathbb{R}^2$ such that
$\partial\, \mathrm{Dsc}(Q)=Q$. Write $D=\mathrm{Dsc}(C)$. Take arcs
$J'$ and $K'$ in $\mathbb{R}^2$ such that
$$\partial J'=\partial K'=\{x,y\}, \hspace{10mm} \mathrm{Int}\,J' \subset \mathrm{Int}~D,$$
$$\mathrm{Int}\,K' \subset \mathbb{R}^2-D, \hspace{5mm}
\mbox{and}\hspace{5mm} \mathrm{Int}\,A' \subset
\mathrm{Int}\,\mathrm{Dsc}(J' \cup K').$$ Then for any $w \in
\mathrm{Int}\,A$, there is a neighborhood $U_w$ of $w$ in
$\mathbb{R}^2$ such that
$$U_w \cap \mathrm{Dsc}(J' \cup K')=\emptyset \hspace{5mm}
\mbox{and} \hspace{5mm} f(U_w) \subset \mathrm{Int}\,\mathrm{Dsc}(J'
\cup K').$$ Let $U=\bigcup\{U_w: w \in \mathrm{Int}\,A\}$. Then $U$
is a neighborhood of $\mathrm{Int}\,A$ in $\mathbb{R}^2$,
$$U \cap \mathrm{Dsc}(J' \cup K')=\emptyset, \hspace{5mm} \mbox{and}
\hspace{5mm} f(U) \subset \mathrm{Int}\,\mathrm{Dsc}(J' \cup K').$$
Clearly, there exist arcs $J$ and $K$ in $U \cup \{x,y\}$ such that
$$\partial J=\partial K=\{x,y\}, \hspace{10mm} \mathrm{Int}\,J \subset \mathrm{Int}\,D,\hspace{10mm} \mathrm{Dsc}(A \cup J) \subset D,$$
$$\mathrm{Int}\,K \subset \mathbb{R}^2-D, \hspace{5mm} \mbox{and} \hspace{5mm} \mathrm{Int}\, A \subset \mathrm{Int}\,\mathrm{Dsc}(J \cup K).$$
Write $E=\mathrm{Dsc}(A \cup K)$. Then
\begin{equation} \label{eq:2-1}
\mathrm{Int}\, E \subset E-A \subset \mathbb{R}^2-D. \tag{2.1}
 \end{equation}
Let the arcs $J''=f(J)$ and $K''=f(K)$. Then $\partial J''=\partial
K''=\{x,y\}$. Since $f$ is an orientation preserving homeomorphism,
we have
$$\mathrm{Int}\, J'' \subset \mathrm{Int}\, D,\hspace{10mm} \mathrm{Int}\, K'' \subset \mathbb{R}^2-D,$$
$$\mathrm{Int}\, A' \subset \mathrm{Int}\, \mathrm{Dsc}(J'' \cup K'').$$
Write $E'=\mathrm{Dsc}(A' \cup K'')$. We have
$$\mathrm{Int}\, E'=f(\mathrm{Int}\, E) \subset f(E-A)=E'-A' \subset \mathbb{R}^2-D,$$
which with \eqref{eq:2-1} implies
\begin{equation} \label{eq:2-2}
f(D) \subset f(\mathbb{R}^2-\mathrm{Int}\,
E)=\mathbb{R}^2-\mathrm{Int}\, E'. \tag{2.2}
 \end{equation}
Write $G=\mathbb{R}^2-\mathrm{Int}\, E'-\mathrm{Int}\,
D-\mathrm{Int}\, A'$. Then $G \cup D=\mathbb{R}^2-\mathrm{Int}\,
E'$, which is closed in $\mathbb{R}^2$. Since $A$ is an absolute
retract, there is a retraction $\gamma: G \rightarrow A$. Since $G
\cap D=A$, we can define a retraction $\beta: G \cup D \rightarrow
D$ by $\beta|G=\gamma$ and $\beta|D=id$. By \eqref{eq:2-2}, we can
define a continuous map $g: D \rightarrow D$ by $g=\beta f$. From
Theorem A we get $\mathrm{Fix}(g) \neq \emptyset$.

\vspace{3mm} For any $w \in \mathrm{Dsc}(A\cup J)$, we have $f(w)
\in \mathrm{Dsc}(A' \cup J'') \subset D$ and $f(w) \neq w$, which
imply $g(w)=f(w) \neq w$.

\vspace{3mm} For any $w \in D-\mathrm{Dsc}(A \cup J)$, if $f(w)
\notin D$ then $g(w)=\beta f(w) \in \beta(G)=\gamma(G)=A \subset
\mathrm{Dsc}(A \cup J)$, which also implies $g(w) \neq w$.

\vspace{3mm} Thus we obtain the following

\vspace{3mm} \textbf{Claim 1.} \ {\it If $v \in \mathrm{Fix}(g)$,
then $v \notin \mathrm{Dsc}(A \cup J)$, and $f(v) \in D$.}

\vspace{3mm} Since $\beta|D=id$, from Claim 1 we see that, if $v \in
\mathrm{Fix}(g)$ then $v=g(v)=\beta f(v)=f(v)$, which implies that
$\mathrm{Fix}(f) \supset \mathrm{Fix}(g) \neq \emptyset$. Lemma 2.1
is proved.  \hfill$\Box$

\vspace{3mm} Basing on Lemma 2.1, we can give a short proof of
proposition C.

\vspace{3mm}\textbf{Proof of Proposition C.}\ Let $\{x,y\}$ be a
$2$-periodic orbit of $f$. Let $r_0=d(x,y)$. If $B^2(x,r_0) \cap
\mathrm{Fix}(f)\neq \emptyset$, then the proposition holds. In the
following we assume that $B^2(x, r_0) \cap
\mathrm{Fix}(f)=\emptyset$.

\vspace{3mm} \textbf{Claim 2.}\ {\it There exists an orientation
preserving homeomorphism $g: \mathbb{R}^2 \rightarrow \mathbb{R}^2$
and an arc $A$ in $\mathbb{R}^2$ such that
$\mathrm{Fix}(g)=\mathrm{Fix}(f)$, $\partial A=\{x_0, y_0\}$,
$f(x_0)=y_0$, $f(y_0)=x_0$, and $g(\mathrm{Int}\, A) \cap
\mathrm{Int}\, A=\emptyset$.}

\vspace{3mm} \textbf{Proof of Claim 2.}\ Let $c=\min\{r>0:
f\left(B^2(x,r)\right) \cap B^2(x,r) \neq \emptyset\}$. Then $0<c<
r_0$. Let $D=B^2(x,c)$ and $E=f(D)$. Then $y \in \mathrm{Int}\, E$,
$E \cap D=\partial E \cap \partial D \neq \emptyset$, and $f(E)\cap
E=f(\partial E) \cap \partial E \neq \emptyset$. Take a point $v \in
\partial E \cap \partial D$, and let $w=f(v)$. Then $w \in f(\partial D)=\partial
E$ and $w \neq v$. Take arcs $A_1$ and $A_2$ in $E$ such that
$$\partial A_1=\{v,y\}, \hspace{8mm} \partial A_2=\{w,y\}, \hspace{8mm} A_1 \cap \partial E=\{v\},$$
$$A_2 \cap \partial E=\{w\}, \hspace{8mm} \mbox{and} \hspace{8mm} A_1 \cap A_2=\{y\}.$$
Let $A_0=f^{-1}(A_2)$ and $A_3=f(A_1)$. Then
$$\partial A_0=\{x,v\}, \hspace{5mm} A_0 \cap \partial D=\{v\}, \hspace{5mm} \partial
A_3=\{x,w\}, \hspace{5mm} \mbox{and} \hspace{5mm} A_3 \cap
E=\{w\}.$$ Let $x_0$ be the point in $A_3 \cap A_0$ such that $[w,
x_0]_{A_3} \cap A_0=\{x_0\}$. Then $x_0 \in [x,v)_{A_0}$. Let
$y_0=f^{-1}(x_0)$. Then $y_0 \in (v, y]_{A_1}$. Write
$$A_0'=[x_0, v]_{A_0}, \hspace{5mm} A_1'=[v,y_0]_{A_1}, \hspace{5mm} \mbox{and} \hspace{5mm} A_3'=[w,x_0]_{A_3}.$$
Then $A_3'=f(A_1')$, and
\begin{equation} \label{eq:2-3}
A_3' \cap A_0'=A_3' \cap A_0=\{x_0\}, \tag{2.3}
 \end{equation}
\begin{equation} \label{eq:2-4}
A_3' \cap A_1' \subset \left((\mathbb{R}^2-E) \cup \{w\}\right) \cap
\left(\mathrm{Int}\, E \cup \{v\}\right)=\emptyset. \tag{2.4}
 \end{equation}
Take homeomorphisms $h_D: \mathbb{R}^2 \rightarrow \mathbb{R}^2$ and
$h_E: \mathbb{R}^2 \rightarrow \mathbb{R}^2$ such that
$$h_D|(\mathbb{R}^2-\mathrm{Int}\, D)=id,\hspace{10mm} h_D(A_0')=A_0,$$
and
$$h_E|(\mathbb{R}^2-\mathrm{Int}\, E)=id,\hspace{10mm} h_E(A_1)=A_1'.$$
Let $A_2'=h_E(A_2)$. Then $A_2' \subset E$, $A_2' \cap \partial
E=\{w\}$, and
\begin{equation} \label{eq:2-5}
A_2' \cap A_0' \subset \left(\mathrm{Int}\, E \cup \{w\}\right) \cap
\left(\mathrm{Int}\, D \cup \{v\}\right)=\emptyset, \tag{2.5}
 \end{equation}
\begin{equation} \label{eq:2-6}
A_2' \cap A_1'=h_E(A_2 \cap A_1)=h_E(\{y\})=\{y_0\}. \tag{2.6}
 \end{equation}
\begin{figure}[!htb]
\centering
\includegraphics[width=13cm]{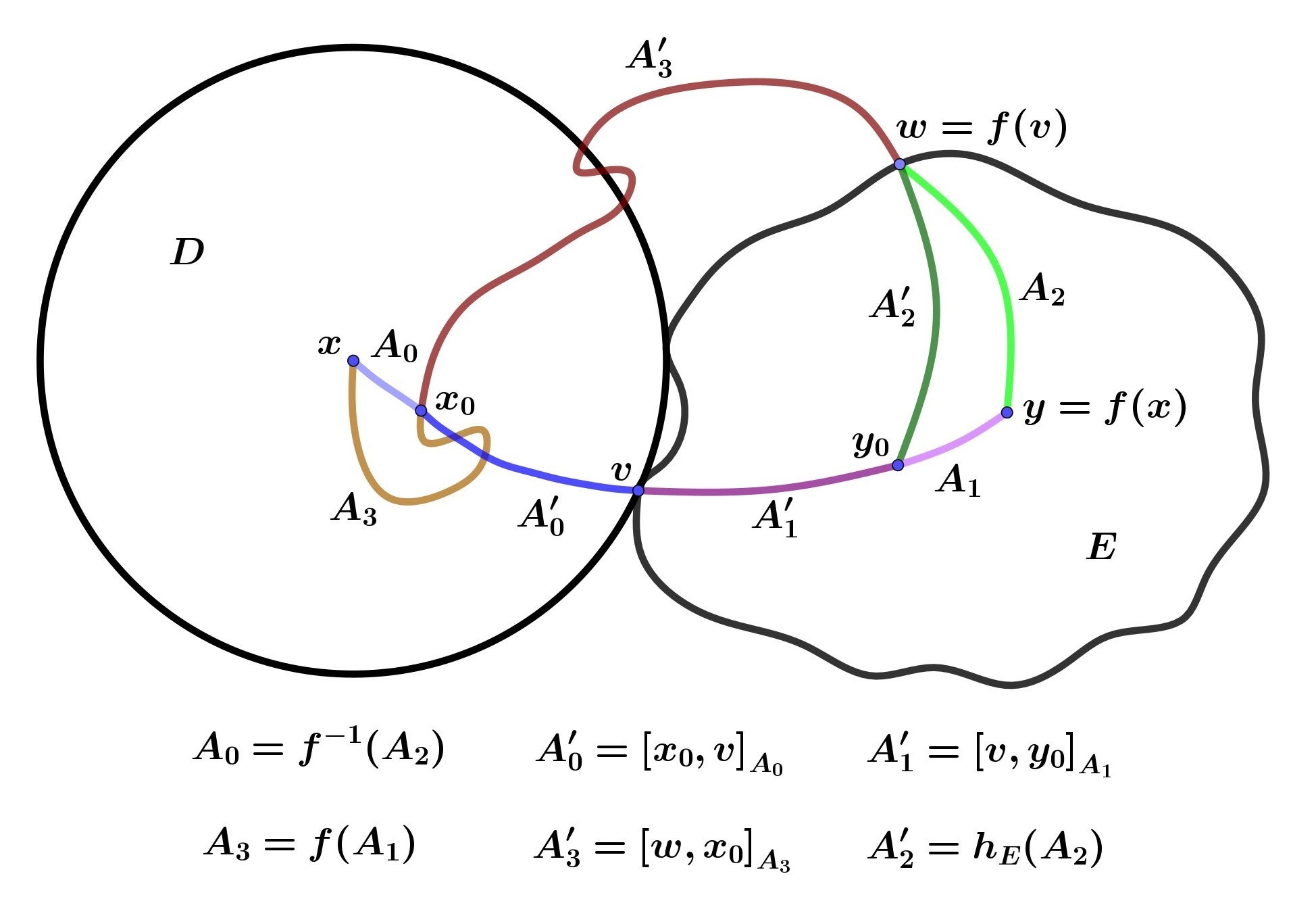}
\caption{}
\end{figure}
\\
Let $g=h_E fh_D: \mathbb{R}^2 \rightarrow \mathbb{R}^2$. Then $g$ is
an orientation preserving homeomorphism, and
\begin{equation} \label{eq:2-7}
g(x_0)=h_Ef(x)=h_E(y)=y_0, \tag{2.7}
 \end{equation}
\begin{equation} \label{eq:2-8}
\begin{split}
& g(y_0)=h_Ef(y_0)=h_E(x_0)=x_0,\\
 & g(A_0')=h_Ef(A_0)=h_E(A_2)=A_2',\\
 & g(A_1')=h_Ef(A_1')=h_E(A_3')=A_3',
 \end{split} \tag{2.8}
 \end{equation}
\begin{equation} \label{eq:2-9}
\mathrm{Fix}(g)=\mathrm{Fix}(f) \subset \mathbb{R}^2-D-E. \tag{2.9}
 \end{equation}
Let $A=A_0' \cup A_1'$. Then $\partial A=\{x_0, y_0\}$,
$$g(A)=g(A_0') \cup g(A_1')=A_2' \cup A_3',$$
and from \eqref{eq:2-3}--\eqref{eq:2-8} we get
$$g(A) \cap A=(A_2' \cup A_3') \cap (A_0' \cup A_1')=\partial A=g(\partial A) \cap \partial A.$$
This means that $g(\mathrm{Int}\, A) \cap \mathrm{Int}\,
A=\emptyset$. Claim 2 is proved.

\vspace{3mm} By Claim 2 and Lemma 2.1 we have $\mathrm{Fix}(g) \neq
\emptyset$, which with \eqref{eq:2-9} implies $\mathrm{Fix}(f) \neq
\emptyset$. Proposition C is proved. \hfill$\Box$

\vspace{3mm} We now give a short proof of Proposition D.

\vspace{3mm}\textbf{Proof of Proposition D.}\ Similar to the proof
of Proposition C, there exists a disc $D=B^2(x,c)$ such that $f(D)
\cap D=f(\partial D) \cap \partial D \neq \emptyset$. For $i \in
\mathbb{N}_n \cup \{0\}$, write $D_i=f^i(D)$ and
$U_i=\,\stackrel{\circ}D_i$. Then $U_i=f^i(\stackrel{\circ}D)$ and
$\partial D_i=f^i(\partial D)$. Take a bounded open set $G$ in
$\mathbb{R}^2$ such that $\bigcup \{D_i: i \in \mathbb{N}_n \cup
\{0\}\} \subset G$. If $G \cap \mathrm{Fix}(f) \neq \emptyset$, then
we can take $g=f$. In the following we assume that $G \cap
\mathrm{Fix}(f)=\emptyset$.

\vspace{3mm} Clearly, there exists a $k \in \mathbb{N}_{n-1}$ such
that $U_{k+1} \cap U_0 \neq \emptyset$ and $U_i \cap U_0=\emptyset$
for $i \in \mathbb{N}_k$. Take a point $v \in U_0$ such that
$f^{k+1}(v) \in U_0$. Write $v_i=f^i(v)$ for $i \in \mathbb{N}_{k+1}
\cup \{0\}$. Then $v_{k+1} \in U_0$, $v_i \in U_i \subset
\mathbb{R}^2-D_0$ for $i \in \mathbb{N}_k$, and $v_i \in U_i \subset
\mathbb{R}^2-D_0-D_1$ for $2 \leq i \leq k$.

\vspace{3mm} If $k=1$, then we can choose an orientation preserving
homeomorphism $h: \mathbb{R}^2 \rightarrow \mathbb{R}^2$ such that
$h(v_2)=v_0$ and $h|(\mathbb{R}^2-U_0)=id$. Put $g=fh$. Then $g$ has
$2$-periodic points $v_1$ and $v_2$, and
$\mathrm{Fix}(g)=\mathrm{Fix}(f)$.

\vspace{3mm} In the following we assume that $2 \leq k \leq n-1$.
Take a point $y \in f(\partial D_0) \cap \partial D_0$. Let $A_0$ be
an arc in $D_0$ and $A_1$ be an arc in $D_1$ such that $\partial
A_0=\{v_{k+1}, y\}$, $\partial A_1=\{y, v_1\}$, and $A_0 \cap
\partial D_0=A_1 \cap \partial D_1=\{y\}$. Let $A=A_0 \cup A_1$. Then
$$A \cap \{f(y), f^{-1}(y)\}=\left\{A-\{y\}\right\} \cap
\left\{f(y), f^{-1}(y)\right\} \subset (U_0 \cup U_1) \cap
\left\{f(y), f^{-1}(y)\right\}=\emptyset,$$ and hence there exists
an open neighborhood $E$ of $y$ in $G-\{v_1, \ldots, v_k\}$ such
that
\begin{equation} \label{eq:2-10}
(E \cup A) \cap \left(f^{-1}(E) \cup f(E)\right)=\emptyset.
\tag{2.10}
 \end{equation}
Noting $A \subset U_0 \cup E \cup U_1$, by \eqref{eq:2-10}, we can
take a disc
\begin{equation} \label{eq:2-11}
K \subset (U_0 \cup E \cup U_1)-\left(f^{-1}(E) \cup \{f(y)\}\right)
\tag{2.11}
 \end{equation}
such that $A \, \subset \; \stackrel{\circ} K$, and $K \cap D_1$ is
a disc, $K \cap \partial D_1$ is an arc. Take an orientation
preserving homeomorphism $h: \mathbb{R}^2 \rightarrow \mathbb{R}^2$
such that
$$h(v_{k+1})=v_1,\hspace{10mm} h(K \cap U_1) \subset K \cap U_1,
\hspace{5mm}\mbox{and} \hspace{5mm} h|(\mathbb{R}^2-K)=id.$$ Let
$g=fh$. Then $g|(\mathbb{R}^2-K)=f|(\mathbb{R}^2-K)$, and hence
\begin{equation} \label{eq:2-12}
\mathrm{Fix}(g) \cap (\mathbb{R}^2-K)=\mathrm{Fix}(f) \cap
(\mathbb{R}^2-K). \tag{2.12}
 \end{equation}

\vspace{3mm} We now consider $\mathrm{Fix}(g) \cap K$. Since $K \cap
f^{-1}(E)=\emptyset$, we have
$$g(K) \cap E=fh(K) \cap E=f(K) \cap E=\emptyset,$$
and hence $K \cap E \cap \mathrm{Fix}(g)=\emptyset$. Since
$$g(K \cap U_1)=fh(K \cap U_1) \subset f(K \cap U_1)
\subset U_2 \subset \mathbb{R}^2-U_1,$$ we have $K \cap U_1 \cap
\mathrm{Fix}(g)=\emptyset$. For any $u \in K \cap U_0$, if $h(u) \in
U_0 \cup U_1$, then
$$g(u) \in f(U_0 \cup U_1) \subset U_1 \cup U_2 \subset \mathbb{R}^2-U_0;$$
if $h(u) \notin U_0 \cup U_1$, then
$$h(u) \in h(K)-(U_0 \cup U_1)=K-(U_0 \cup U_1) \subset E,$$
which with \eqref{eq:2-11} implies $g(u)=fh(u) \in f(E) \subset
\mathbb{R}^2-K$. Thus $K \cap U_0 \cap \mathrm{Fix}(g) =\emptyset$.
To sum up, we get $K \cap \mathrm{Fix}(g)=\emptyset$, which with
\eqref{eq:2-12} and $K \subset G \subset
\mathbb{R}^2-\mathrm{Fix}(f)$ implies
$\mathrm{Fix}(g)=\mathrm{Fix}(f)$.

\vspace{3mm} In addition, since $K \cap \{v_2, \ldots,
v_k\}=\emptyset$, we have $g(v_{k+1})=fh(v_{k+1})=f(v_1)=v_2$ and
$g(v_i)=f(v_i)=v_{i+1}$ for $i \in \{2, \ldots, k\}$. Thus $g$ has
$k$-periodic points. Proposition D is proved. \hfill$\Box$

\section{Orientation preserving homeomorphisms
of the plane having BP-chain recurrent points}

Let $X$ be a metric space, $f \in C^0(X)$, and $x \in X$.

\vspace{3mm} $x$ is called an {\it \textbf{almost periodic point}}
of $f$ if for any neighborhood $U$ of $x$ there exists $N \in
\mathbb{N}$ such that for any $k \in \mathbb{N}$, one has
$\{f^{k+i}(x): i=0, 1, \ldots, N\} \cap U \neq \emptyset$.

\vspace{3mm} $x$ is called a {\it \textbf{recurrent point}} of $f$,
if for any neighborhood $U$ of $x$ in $X$, there exists $n \in
\mathbb{N}$ such that $f^n(x) \in U$.

\vspace{3mm} $x$ is called an {\it \textbf{$\omega$-limit point}} of
$f$ if there exist $y \in X$ and positive integers
$n_1<n_2<n_3<\cdots$ such that $\lim_{k \rightarrow
\infty}f^{n_k}(y)=x$.

\vspace{3mm} $x$ is called a {\it \textbf{non-wandering point}} of
$f$ if for any neighborhood $U$ of $x$ in $X$, there exists $n \in
\mathbb{N}$ such that $f^n(U) \cap U \neq \emptyset$.

\vspace{3mm} For any $\varepsilon>0$, $n \in \mathbb{N}$ and $x, y
\in X$, a sequence $(x_0, x_1, \ldots, x_n)$ of points in $X$ is
called an {\it \textbf{$\varepsilon$-chain}} of $f$ from $x$ to $y$
if $x_0=x$, $x_n=y$, and $d(f(x_{k-1}),x_k)<\varepsilon$ for any $k
\in \mathbb{N}_n$.

\vspace{3mm} $x \in X$ is called a {\it \textbf{chain recurrent
point}} of $f$ if for any $\varepsilon>0$, there exists an
$\varepsilon$-chain of $f$ from $x$ to $x$.

\vspace{3mm} Denote by $P(f)$, $AP(f)$, $R(f)$, $\omega(f)$,
$\Omega(f)$ and $CR(f)$ the sets of periodic points, almost periodic
points, recurrent points, $\omega$-limit points, non-wandering
points and chain recurrent points of $f$, respectively. Obviously,
we have $$\mathrm{Fix}(f) \subset P(f) \subset AP(f) \subset R(f)
\subset \omega(f) \subset \Omega(f) \subset CR(f).$$

In \cite{Fa}\;\!, J. Fathi proved the following theorem. This means
that the condition ``having a periodic point'' in Theorem B can be
weakened to ``having a non-wandering point''.

\vspace{3mm}\textbf{Theorem E.}\ {\it Let $f:\Bbb R^2\to \Bbb R^2$
be an orientation preserving homeomorphism. If $f$ has a
non-wandering point, then $f$ has a fixed point}\;\!.

\vspace{3mm} We now give an example of orientation preserving
homeomorphisms $f: \mathbb{R}^2 \rightarrow \mathbb{R}^2$, which has
no non-wandering points, but the all points in $\mathbb{R}^2$ are
chain recurrent points of $f$.

\vspace{3mm}\textbf{Example 3.1.}\ For convenience, we regard the
real line $\mathbb{R}$ as a subspace of the plane $\mathbb{R}^2$,
that is, any point $r \in \mathbb{R}$ and the point $(r,0) \in
\mathbb{R}^2$ are regarded as the same.

\vspace{3mm} Define projections $p: \mathbb{R}^2 \rightarrow
\mathbb{R}$ and $q: \mathbb{R}^2 \rightarrow \mathbb{R}$ and a
reflection $\varphi: \mathbb{R}^2 \rightarrow \mathbb{R}^2$ by, for
any $(r,s) \in \mathbb{R}^2$,
$$p(r,s)=r,\hspace{10mm} q(r,s)=s, \hspace{10mm} \mbox{and} \hspace{10mm} \varphi(r,s)=(r,-s).$$
For any $t \in (0, +\infty)$, let
$$K_{t1}=\left\{(r,s) \in \mathbb{R}^2: r \in (0, +\infty), s=t/r\right\},
\hspace{10mm} K_{t2}=\varphi(K_{t1}).$$ Then $K_{t1}$ and $K_{t2}$
are hyperbolas. For any $t \in [1, +\infty)$, let
$$C_t=S^1\left(0, \sqrt{2t}\right), \hspace{10mm} w_{t1}=(\sqrt{t},\sqrt{t}),\hspace{10mm} w_{t2}=\varphi(w_{t1}).$$
For any $t \in (0,1]$, let
$$\lambda(t)=(1-t^6)/t,\hspace{10mm} \mu(t)=t^2\sqrt{1+t^6},$$
and let
$$C_t=S^1(\lambda(t), \mu(t)), \hspace{10mm} w_{t1}=\left(1/t, t^2\right),
\hspace{10mm} w_{t2}=\varphi(w_{t1}).$$ Note that the center, of the
circle $C_t$, $\lambda_t \in \mathbb{R}=\mathbb{R} \times \{0\}
\subset \mathbb{R}^2$. It is easy to check that, for each $t \in (0,
+\infty)$,
$$C_t \cap K_{t1}=\{w_{t1}\}, \hspace{10mm} C_t \cap K_{t2}=\{w_{t2}\},$$
and the tangent lines of $C_t$ and $K_{t1}$ at the point $w_{t1}$
are the same, the tangent lines of $C_t$ and $K_{t2}$ at the point
$w_{t2}$ are the same.

\vspace{3mm} For each $t \in (0, +\infty)$, write
$$A_t=\{z \in C_t: p(z) \leq p(w_{t1})\},$$
$$J_{t1}=\{z \in K_{t1}: p(z) \geq p(w_{t1})\}, \hspace{10mm} J_{t2}=\varphi(J_{t1}),$$
$$L_t=J_{t1} \cup A_t \cup J_{t2}.$$
Then $L_t$ is a smooth curve, $\bigcup \{L_t: t \in (0,
+\infty)\}=\mathbb{R}^2$, and for $0<t<t'<+\infty$, $L_t \cap
L_{t'}=\emptyset$.

\begin{figure}[!htb]
\centering
\includegraphics[width=10cm]{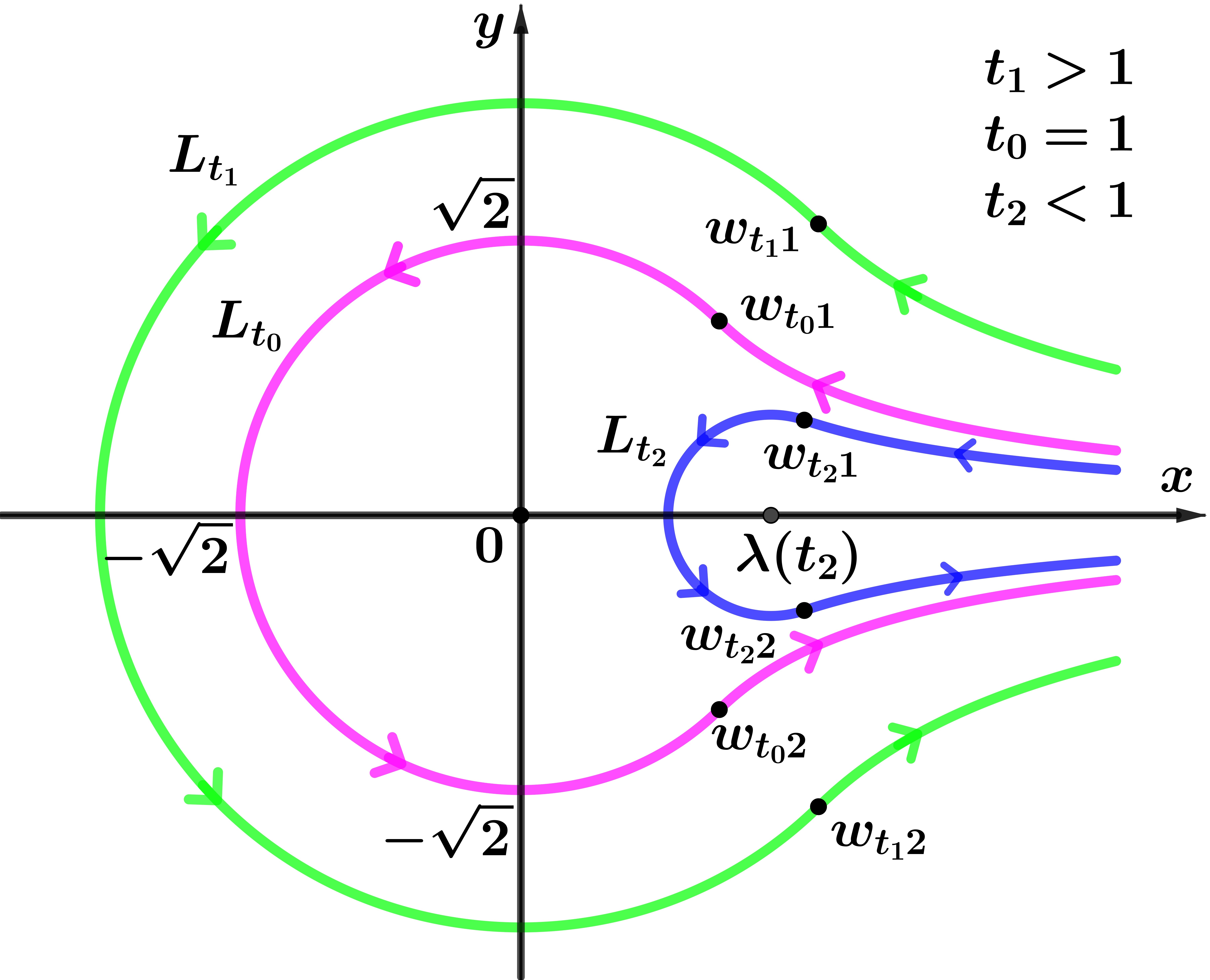}
\caption{}
\end{figure}

\vspace{3mm} Let $v_{t1}$ be the point in $A_t$ such that
$q(v_{t1})>0$ and $p(v_{t1})=p(w_{t1})-q(w_{t1})/2$, and let
$v_{t2}=\varphi(v_{t1})$. Write
$$A_{t1}=A_t[v_{t1}, w_{t1}], \hspace{10mm} A_{t2}=\varphi(A_{t1}), \hspace{10mm}
\mbox{and} \hspace{10mm} A_{t0}=A_t[v_{t1}, v_{t2}].$$ Clearly,
there exists a unique orientation preserving homeomorphism $f:
\mathbb{R}^2 \rightarrow \mathbb{R}^2$ such that, for any $t \in (0,
+\infty)$,

\vspace{3mm} (a)\ $f(L_t)=L_t$, $f(J_{t1})=J_{t1} \cup A_{t1}$, and
$$f(A_{t1} \cup A_{t0})=A_{t0} \cup A_{t2}, \hspace{10mm} f(A_{t2} \cup J_{t2})=J_{t2};$$

(b)\ For any $z \in J_{t1}$, $p(z)-p(f(z))=q(z)/2$. Specially, we
have $f(w_{t1})=v_{t1}$;

\vspace{3mm} (c)\ For any $z \in A_{t2} \cup J_{t2}$, $f(z)=\varphi
f^{-1}\varphi(z)$. Specially, we have $f(v_{t2})=w_{t2}$;

\vspace{3mm} (d)\ For any $z \in A_{t1} \cup A_{t0}$, $d(f(z),
z)=d(v_{t1}, w_{t1})$.

\vspace{3mm} It is easy to see that, in Example 3.1, the
homeomorphism $f: \mathbb{R}^2 \rightarrow \mathbb{R}^2$ has no
non-wandering points, but the all points in $\mathbb{R}^2$ are chain
recurrent points of $f$. From this example we see that the condition
``having a periodic point'' in Theorem B can not be weaken to
``having a chain recurrent point''.

\vspace{3mm} We now propose a notion, which is weaker than the
notion of non-wandering points, but is stronger than that of chain
recurrent points.

\vspace{3mm}\textbf{Definition 3.2.} \ Let $(X,d)$ be a metric
space, and $f: X \rightarrow X$ be a continuous map. A point $x \in
X$ is called a {\it \textbf{chain recurrent point of $f$ with
bounded perturbation}} (for short, a {\it \textbf{BP-chain recurrent
point of $f$}}) if there exists a bounded subset $W$ of $X$ such
that, for any $\varepsilon>0$, there is an $\varepsilon$-chain
$(x_0, x_1, \ldots, x_n)$ of $f$ from $x$ to $x$ such that, for any
$k \in \mathbb{N}_n$, if $f(x_{k-1}) \in X-W$ then $x_k=f(x_{k-1})$.

\vspace{3mm} Note that, in Definition 3.2, the condition ``if
$f(x_{k-1}) \in X-W$ then $x_k=f(x_{k-1})$'' is equivalent to that
``if $x_k \neq f(x_{k-1})$ then $f(x_{k-1}) \in W$''.

\vspace{3mm} Let $y, z \in X$ and $(y_0, y_1, \ldots, y_m)$ be an
$\varepsilon$-chain of $f$ from $y$ to $z$. For $k \in
\mathbb{N}_m$, if $d(f(y_{k-1}), y_k)>0$ (resp., $d(f(y_{k-1}),
y_k)=0$) then we say that the $\varepsilon$-chain $(y_0, y_1,
\ldots, y_m)$ has a {\it \textbf{nonzero perturbation}} (resp., has
a {\it \textbf{zero perturbation}}, or has no {\it
\textbf{perturbation}}) at the point $y_k$. By this notion, the
$\varepsilon$-chain $(x_0, x_1, \ldots, x_n)$ in Definition 3.2 has
no perturbation at $x_k$ if $f(x_{k-1}) \notin W$.

\vspace{3mm} Denote by $CR_{BP}(f)$ the set of all BP-chain
recurrent points of $f$. Obviously, we have $CR_{BP}(f) \subset
CR(f)$, and if $X$ is a bounded metric space then
$CR_{BP}(f)=CR(f)$.

\vspace{3mm}\textbf{Proposition 3.3.} \ {\it Let $(X, d)$ be a
metric space, and $f: X \to X$ be a continuous map. Then $\Omega(f)
\subset CR_{BP}(f)$}\;\!.

\vspace{3mm}\textbf{Proof.} \ Consider any $x \in \Omega(f)$. If $x
\in \mathrm{Fix}(f)$ then it is clear that $x \in CR_{BP}(f)$. We
now assume that $x \notin \mathrm{Fix}(f)$. Let $c=d(x, f(x))$, and
take $W=B(x, 2c)$. For any given $\varepsilon \in (0, c]$, take a
$\delta \in (0, \varepsilon/3]$ such that
$$f(x, B(x, \delta)) \subset B(f(x), \varepsilon/3).$$
Since $x \in \Omega(f)$, there exist $y \in B(x, \delta)$ and $n
\geq 2$ such that $f^n(y) \in B(x, \delta)$. Let $x_0=x_n=x$, and
let $x_k=f^k(y)$ for $k \in \mathbb{N}_{n-1}$. Then $(x_0, x_1,
\ldots, x_n)$ is an $\varepsilon$-chain of $f$ from $x$ to $x$ which
has no perturbation at $x_k$ if $f(x_{k-1}) \in X-W$. \hfill$\Box$

\vspace{3mm} It is easy to see that, in Example 3.1,
$CR_{BP}(f)=\emptyset$, although $CR(f)=\mathbb{R}^2$. In the
following example, we choose an unbounded open subset $U$ of
$\mathbb{R}^2$ and give a homeomorphism $f: U \rightarrow U$ such
that $\Omega(f)=\emptyset$ and $CR_{BP}(f)=U$.

\vspace{3mm} \textbf{Example 3.4.}\ Let $U=\mathbb{R}^2-\mathbb{R}$.
Write $x=-1$, $y=1$. Note that $\{x,y\} \subset [-1,1] \subset
\mathbb{R}=\mathbb{R} \times \{0\} \subset \mathbb{R}^2$. For any $t
\in \mathbb{R}$, let the point $w_t=(0,t) \in \mathbb{R}^2$, let
$r_t=d(w_t, y)=\sqrt{t^2+1}$, and let the circle $C_t=S^1(w_t,
r_t)$. Then $\{x,y\} \subset C_t$. Let $q: \mathbb{R}^2 \rightarrow
\mathbb{R}$ and $\varphi: \mathbb{R}^2 \rightarrow \mathbb{R}^2$ be
the same as in Example 3.1. Write
$$A_{t1}=\{z \in C_t: q(z) \geq 0\}, \hspace{10mm} A_{t2}=\{z \in C_t: q(z) \leq 0\}.$$
Then $A_{t1}$ and $A_{t2}$ are arcs,
$$\partial A_{t1}=\partial A_{t2}=\{x, y\}, \hspace{10mm} A_{t2}=\varphi(A_{-t,1}),$$
$$\bigcup \{\mathrm{Int}\, A_{t1} \cup \mathrm{Int}\, A_{t2}: t \in \mathbb{R}\}=U,$$
and for $-\infty<t<t'<+\infty$, the open arcs $\mathrm{Int}\,
A_{t1}$, $\mathrm{Int}\, A_{t2}$, $\mathrm{Int}\, A_{t'1}$ and
$\mathrm{Int}\, A_{t'2}$ are pairwise disjoint. Clearly, there
exists a unique homeomorphism $f: U \rightarrow U$ such that, for
any $t \in \mathbb{R}$,

\vspace{3mm} (a)\ $f(\mathrm{Int}\, A_{t1})=\mathrm{Int}\, A_{t1}$,
and $f(\mathrm{Int}\, A_{t2})=\mathrm{Int}\, A_{t2}$;

\vspace{3mm} (b)\ For any $z \in A_{t1} \cup A_{t2}$,
$d(f(z),z)=d(z, \{x, y\})/2$, and if $z \in A_{t1}$ then $f(z) \in
A_{t1}(z,x)$, if $z \in A_{t2}$ then $f(z) \in A_{t2}(z,y)$.

\vspace{3mm} It is easy to see that, in Example 3.4, we have
$\Omega(f)=\emptyset$ and $CR_{BP}(f)=U$.

\vspace{3mm} The condition ``having a periodic point'' in Theorem B
and the condition ``having a non-wandering point'' in Theorem E can
be weakened to ``having a BP-chain recurrent point''. In fact, we
have the following

\vspace{3mm}\textbf{Theorem 3.5.} \ {\it Let $f:\mathbb{R}^2\to
\mathbb{R}^2$ be an orientation preserving homeomorphism. If $f$ has
a BP-chain recurrent point, then $f$ has a fixed point}\;\!.

\vspace{3mm}\textbf{Proof.} \ Let $x$ be a BP-chain recurrent point
of $f$, and let the bounded set $W \subset \mathbb{R}^2$ be the same
as in Definition 3.2. We may assume that $x \in W$. Take discs $D'$
and  $D$ in $\mathbb{R}^2$ such that
$$f^{-1}(W) \cup  W \cup f(W) \subset D' \subset B^2(D', 1) \subset D.$$
Let
$$c=\min\{d(f(y),y), d(f^{-1}(y), y): y \in D\}.$$
If $\mathrm{Fix}(f)=\emptyset$, then
$\mathrm{Fix}(f^{-1})=\emptyset$, and hence we have $c>0$. Take a
$\delta_0 \in (0, c/3]$ such that, for any $y \in D$,
$$f(B^2(y, \delta_0)) \subset B^2\left(f(y), c/3\right).$$
Take an $\varepsilon \in (0, \delta_0/2]\cap (0,1]$ such that, for
any $y \in D$ and $i \in \{-1,1\}$,
\begin{equation} \label{eq:3-1}
f^i\left(B^2(y, 2\varepsilon)\right) \subset B^2\left(f^i(y),
\delta_0\right). \tag{3.1}
\end{equation}
By Definition 3.2, there is an $\varepsilon$-chain $(x_0, x_1,
\ldots, x_n)$ of $f$ from $x$ to $x$ such that, for any $k \in
\mathbb{N}_n$, if $f(x_{k-1}) \neq x_k$ then $f(x_{k-1}) \in W$.

\vspace{3mm} In this $\varepsilon$-chain, if $x_j=x_k$ for some
$0<j<k \leq n$, then we we can replace $(x_0, x_1, \ldots, x_n)$ by
$(x_0, \ldots, x_{j-1}, x_k, \ldots, x_n)$. Thus we may assume $x_j
\neq x_k$ for any $0<j<k \leq n$. Let
$$L_0=\{x\}, \hspace{5mm} \mbox{and} \hspace{5mm} L_k=[f(x_{k-1}), x_k] \mbox{\ for\ } k \in \mathbb{N}_n.$$
Then we have

\vspace{3mm} (P.1)\, For $k \in \mathbb{N}_n$, $L_k \neq \{x_k\}$ if
and only if $f(x_{k-1}) \neq x_k$;

\vspace{3mm} (P.2)\, $0 \leq \mathrm{diam}\, L_k<\varepsilon \leq
\delta_0/2 \leq c/6$, for any $k \in \mathbb{N}_n \cup \{0\}$.

\vspace{3mm} (P.3)\, For $k \in \mathbb{N}_n \cup \{0\}$, if $L_k
\neq \{x_k\}$, then $k \geq 1$, $f(x_{k-1}) \in W \subset D$, and
for any disc $E_k$ in $\mathbb{R}^2$ with $L_k \subset \,
\stackrel{\circ} E_k$ and $\mathrm{diam}\, E_k<2\varepsilon$, from
\eqref{eq:3-1} we get
$$d\left(f^i(E_k), E_k\right)>d\left(f^i(f(x_{k-1})), f(x_{k-1})\right)-2\varepsilon-\delta_0 \geq c-2\delta_0 \geq c/3>0$$
for $i \in \{-1,1\}$. In addition, if $L_k=\{x_k\}$, then there is
also a disc $E_k$ in $\mathbb{R}^2$ such that $x_k \in \,
\stackrel{\circ}E_k$ and $d(f^i(E_k), E_k)>0$ for $i \in \{-1,1\}$.

\vspace{3mm} \textbf{Claim 3.}\ {\it $d(L_{k-1}, L_k)>0$ for any $k
\in \mathbb{N}_n$.}

\vspace{3mm} \textbf{Proof of Claim 3.}\  {\it\textbf{Case $1$.}}\,
If $L_{k-1}=\{x_{k-1}\}$ and $L_k=\{x_k\}$, then $d(L_{k-1}, L_k)>0$
since $x_k \neq x_{k-1}$.

\vspace{3mm} {\it\textbf{Case $2$.}}\, If $L_{k-1}=\{x_{k-1}\}$ and
$L_k \neq \{x_k\}$, then $f(x_{k-1}) \neq x_k$, $f(x_{k-1}) \in W$,
$d(x_{k-1}, f(x_{k-1})) \geq c$, and hence
$$d(L_{k-1}, L_k)=d(x_{k-1}, [f(x_{k-1}), x_k]) \geq c-\mathrm{diam}\,L_k>c-\varepsilon>2c/3>0.$$

\vspace{3mm} {\it\textbf{Case $3$.}}\, If $L_{k-1} \neq
\{x_{k-1}\}$, then, no matter whether $L_k=\{x_k\}$ or not, we have
$$k \geq 2, \hspace{10mm} f(x_{k-2}) \neq x_{k-1}, \hspace{5mm} \mbox{and} \hspace{5mm} f(x_{k-2}) \in W.$$
Write
$$x_{k-1}'=f(x_{k-2}), \hspace{10mm} x_k'=f(x_{k-1}), \hspace{10mm} x_k''=f(x_{k-1}').$$
Then $$L_{k-1}=[x_{k-1}', x_{k-1}], \hspace{10mm} L_k=[x_k', x_k],$$
$$x_{k-1}' \in D, \hspace{10mm} d(x_{k-1}', x_k'') \geq c,$$
$$d(x_{k-1}', x_{k-1})<\varepsilon \leq \delta_0/2, \hspace{10mm} d(x_k'', x_k') \leq c/3,$$
$$d(x_{k-1}', x_k') \geq c-c/3=2c/3,$$
and hence, in Case 3 we also have
$$d(L_{k-1}, L_k) \geq d(x_{k-1}', x_k')-\mathrm{diam}\, L_{k-1}-\mathrm{diam}\, L_k>2c/3-2\varepsilon \geq c/3>0.$$
Claim 3 is proved.

\vspace{3mm} Since $x_n=x_0=x$, we have $L_n \cap L_0 \neq
\emptyset$. Hence, by (P.2), (P.3) and Claim 3, there exist integers
$p$ and $q$ such that

\vspace{3mm} (P.4)\, $0 \leq p<p+1<q\leq n$, $L_q \cap L_p \neq
\emptyset$, $\mathrm{diam}(L_q \cup L_p)<2\varepsilon$, and $L_k
\cap L_j=\emptyset$ for any integers $\{j, k\} \subset [p,q]$ with
$0<k-j<q-p$;

\vspace{3mm} (P.5)\, There exist pairwise disjoint disc $D_{p+1},
D_{p+2}, \ldots, D_q$ such that

\vspace{3mm} (a)\, $L_k \subset\, \stackrel{\circ}D_k$ and
$\mathrm{diam}\, D_k<\varepsilon$ for $p<k<q$;

\vspace{3mm} (b)\, $L_p \cup L_q \subset\, \stackrel{\circ}D_q$, and
$\mathrm{diam}\, D_q<2\varepsilon$;

\vspace{3mm} (c)\, For $p<k \leq q$, $d(f^{-1}(D_k), D_k)>0$.
Particularly, if $D_k \cap D \neq \emptyset$ then $$d(f^{-1}(D_k),
D_k)>c/4.$$

For $p<k \leq q$, let $G_{k-1}=f^{-1}(D_k)$. Then $G_{k-1} \cap
D_k=\emptyset$, $x_{k-1} \in f^{-1}(L_k) \subset \, \stackrel{\circ}
G_{k-1}$, and $G_{p+1}, G_{p+2}, \ldots, G_q$ are pairwise disjoint
discs. Take a homeomorphism $h_{k-1}: G_{k-1} \rightarrow D_k$ such
that

\vspace{3mm} (i)\, $h_{k-1}|\partial G_{k-1}=f|\partial G_{k-1}$;

\vspace{3mm} (ii)\, $h_{k-1}(x_{k-1})=x_k$ for $p<k<q$, and
$h_{q-1}(x_{q-1})=x_p$.

\vspace{3mm} Let $X=\mathbb{R}^2-\bigcup\{G_{k-1}: p<k \leq q\}$.
Define a map $h: \mathbb{R}^2 \rightarrow \mathbb{R}^2$ by
$h|X=f|X$, and $h|{G_{k-1}}=h_{k-1}$ for $p<k\leq q$. Then $h$ is an
orientation preserving homeomorphism, and
$\mathrm{Fix}(h|X)=\mathrm{Fix}(f|X)=\emptyset$. For $p<k \leq q$,
from $G_{k-1}\cap D_k=\emptyset$ we get
$\mathrm{Fix}(h|G_{k-1})=\mathrm{Fix}(h_{k-1})=\emptyset$. Thus
$\mathrm{Fix}(h)=\mathrm{Fix}(f)=\emptyset$. On the other hand, $h$
has a periodic orbit $\{x_p, x_{p+1}, \ldots, x_{q-1}\}$. By Theorem
B, we have $\mathrm{Fix}(h) \neq \emptyset$. This leads to a
contradiction. Thus we must have $\mathrm{Fix}(f) \neq \emptyset$.
Theorem 3.5 is proved. \hfill$\Box$

\subsection*{Acknowledgements}
The authors are supported by NNSF of China (Grant No. 12261006) and
Project of Guangxi First Class Disciplines of Statistics (No.
GJKY-2022-01). The second author is also supported by NNSF of China
(Grant No. 12171175). We thank Hui Xu for his help in drawing Figure
2.


\end{document}